\documentstyle[12pt]{article} 
\begin{document} 
\newtheorem{Def}{Definition}[section] 
\newtheorem{thm}{Theorem}[section] 
\newtheorem{lem}{Lemma}[section] 
\newtheorem{rem}{Remark}[section] 
\newtheorem{prop}{Proposition}[section] 
\newtheorem{cor}{Corollary}[section] 
\newtheorem{conj}{Conjecture}[section]
\newtheorem{question}{Question}[section]
\title 
{Degenerate conformally invariant
fully nonlinear elliptic equations }
\author{YanYan Li\thanks{Partially
 supported by
       NSF grant DMS-0401118.}\\
       Department of Mathematics\\
       Rutgers University\\
       110 Frelinghuysen Road\\
       Piscataway, NJ 08854\\
       USA
       }
\date{} 
\maketitle 
\input { amssym.def} 

\setcounter{section}{0}

\section{Introduction}

There has been much work on conformally invariant
fully nonlinear elliptic equations
and  applications to geometry and topology.
See for instance \cite{V1},
 \cite{CGY}, \cite{C}, \cite{LL1},  \cite{LL3},
\cite{GV}, and the references therein.
An important issue in the study of such equations
is to  classify entire solutions which
arise from  rescaling blowing up solutions.  Liouville type
theorems for general conformally invariant 
fully nonlinear second order elliptic equations 
have been obtained in \cite{LL3}. For previous works on the
subject, see  \cite{LL3} for a description.
Classification of entire solutions 
to degenerate equations is also of importance, as demonstrated in
 \cite{CGY2}.  In this paper we 
give  Liouville type theorems for general
degenerate  conformally invariant
fully nonlinear  second order elliptic equations.

Let ${\cal S}^{n\times n}$ denote the set of
$n\times n$ real symmetric matrices,
 ${\cal S}^{n\times n}_+$ denote the subset of ${\cal S}^{n\times n}$
 consisting of positive definite matrices,
  $O(n)$ denote
  the set of $n\times n$ real orthogonal matrices,
   $U\subset {\cal S}^{n\times n}$
   be  an open set, and  $F\in C^1(U)\cap C^0(\overline U)$.
   
   We list below a number of properties of $(F, U)$.
   Subsets of these properties
   will be used in various lemmas, propositions  and theorems:  
\begin{equation}
O^{-1}UO=U,\qquad
\forall\ O\in O(n),
\label{21-1}
\end{equation}

\begin{equation}
M\in U\ \mbox{and}\ N\in {\cal S}^{ n\times n}_+\
\mbox{implies}\ M+N\in U,
\label{C3-1}
\end{equation}

\begin{equation}
M\in U\ \mbox{implies}\ aM\in U\
\mbox{for all positive constant}\ a,
\label{A1-0}
\end{equation}

\begin{equation}
\{aI\ |\ a>0\}\cap \partial U=\emptyset,
\label{A3-0}
\end{equation}
where $I$ denotes the $n\times n$ identity matrix.

Let $F\in C^1(U)\cap C^0(\overline U)$ satisfy
\begin{equation}
F(O^{-1}MO)=F(M),\quad
\forall\ M\in U, \forall\ O\in O(n),
\label{21-3}
\end{equation}
\begin{equation}
(F_{ij}(M)))>0,\qquad
\forall\ M\in U,
\label{21-4}
\end{equation}
where $F_{ij}(M):=\frac{\partial F}{ \partial M_{ij} }(M)$, and,
\begin{equation}
F>0\ \mbox{in}\ U,\qquad
F=0\ \mbox{on}\ \partial U,
\label{A1-1}
\end{equation}

Examples of such $(F,U)$ include those given by the
elementary symmetric functions. 
 For  $1\le k\le n$, let
 $$
 \sigma_k(\lambda)=\sum_{1\le i_1<\cdots <i_k\le n}\lambda_{i_1}
 \cdots \lambda_{i_k}
 $$
 be the $k-$th elementary symmetric function and let
 $
 \Gamma_k$ be the connected component
 of $\{\lambda\in \Bbb R^n\ |\
 \sigma_k(\lambda)>0\}$ containing the positive cone $\Gamma_n
 :=\{\lambda=(\lambda_1, \cdots, \lambda_n)\ |\
 \lambda_i>0\}$.
Let
$$
U_k:=\{ M\in {\cal S}^{n\times n}\
|\ \lambda(M)\in \Gamma_k\},
$$ 
and
$$
F_k(M):=\sigma_k(\lambda(M))^{ \frac 1k},
$$
where $\lambda(M)$ denotes the eigenvalues of $M$.
Then $(F,U)=(F_k, U_k)$ satisfy all the above listed 
properties, see for instance \cite{CNS}.

Other, much more general, examples are as follows.
Let
$$
\Gamma\subset \Bbb R^n \ \mbox{be
an open convex symmetric cone  with vertex at the origin }
$$
satisfying
$$
\Gamma_n
\subset \Gamma\subset \Gamma_1:=\{\lambda\in\Bbb{R}^n|\sum\limits_i
\lambda_i>0\}.
$$
Naturally, $\Gamma$ being symmetric
means $(\lambda_1, \lambda_2,\cdots, \lambda_n)\in
\Gamma$ implies $(\lambda_{i_1}, \lambda_{i_2}, \cdots , \lambda_{i_n})
\in \Gamma$ for any permutation
 $(i_1, i_2, \cdots , i_n)$ of
 $(1,2,\cdots, n)$.

Let
$$
f\in C^1(\Gamma)\cap
 C^0(\overline \Gamma)
 $$
 satisfy
 $$
 f|_{\partial\Gamma}=0,\qquad
 \nabla f\in\Gamma_n\ \mbox{on}\ \Gamma,
 $$
 and
$$
 f(s\lambda)=sf(\lambda),\qquad \forall\ s>0\
 \mbox{and}\ 
  \lambda\in \Gamma.
 $$

With such $(f, \Gamma)$,  let
$$
U:=\{M\in {\cal S}^{n\times n}\ |\
\lambda(M)\in \Gamma\},
$$
and
$$
F(M):=f(\lambda(M)).
$$
Then $(F, U)$ satisfies all the above listed properties.
In fact, for all these $(F,U)$, $A^u\in U$ implies
$\Delta u\le 0$ --- see below for the definition of $A^u$.  
So for these $(F,U)$, the assumption
$\Delta u\le 0$ in various theorems in this paper is automatically
satisfied.
We note that in all these examples, $F$ is actually concave in $U$,
but this property is not needed for results in this paper.

As mentioned above, entire solutions to general equation
$$
F\left(A^u\right)=1, \ u>0,\ A^u\in U, \qquad \mbox{in}\ \Bbb R^n
$$
are classified in 
 \cite{LL3}.
 Here and throughout the paper we use notation
 $$
 A^u= -\frac{2}{n-2}u^{  -\frac {n+2}{n-2} }
 \nabla^2u+ \frac{2n}{(n-2)^2}u^ { -\frac {2n}{n-2} }
 \nabla u\otimes\nabla u-\frac{2}{(n-2)^2} u^ { -\frac {2n}{n-2} }
 |\nabla u|^2I,
 $$
where $\nabla u$ denotes the gradient of $u$ and
$\nabla^2 u$ denotes the Hessian of $u$.

In this paper we classify appropriate weak solutions to
$$
F\left(A^u\right)=0, \ u>0, \qquad \mbox{in}\ \Bbb R^n.
$$
The techniques developed in \cite{LL3} play 
important roles in our studies.
As in  \cite{LL3}, we make use of the method of moving
spheres,
a variant of the method of moving
planes which fully exploits the conformal invariance of the
problem.
The method of moving planes has been used in
classical works of Gidas, Ni and Nirenberg
\cite{GNN} and Caffarelli, Gidas and Spruck
\cite{CGS}, and others --- see for instance
\cite{LL3} for a description, to study Liouville type theorems.

For $x\in \Bbb R^n$, $\lambda>0$, and for some function
$u$, we denote the Kelvin transformation of $u$ with
respect to $B_\lambda(x)$ by
$$
 u_{x,\lambda}(y):=\frac { \lambda^{n-2} }
 { |y-x|^{n-2} }  u(x+ \frac { \lambda^2 (y-x) }{ |y-x|^2 } ).
 $$
Here and throughout
 the paper we use $B_a(x)\subset \Bbb R^n$ to denote 
the ball of radius $a$ and centered at  $x$, and use
$B_a$ to denote $B_a(0)$. 
Also, unless otherwise stated, the dimension $n$ is bigger than $2$.

We first introduce a notion of weak solutions to the degenerate
equations.

\begin{Def} Let $U\subset {\cal S}^{n\times n}$ be an open
set and  $F\in C^1(U)\cap C^0(\overline U)$
 satisfy 
(\ref{21-1}),
 (\ref{C3-1}), 
(\ref{A1-0}), 
(\ref{21-3}), 
(\ref{21-4}), 
and (\ref{A1-1}).
A positive continuous function
$u$ on an open set $\Omega$ of $\Bbb R^n$ is said to be a weak solution of
\begin{equation}
F(A^u)=0\qquad
\mbox{in}\ \Omega
\label{A1-2}
\end{equation}
if there exist 
$\{u_i\}$ in $C^2(\Omega)$  and $\{\beta_i\}$ in $C^0(\Omega,
{\cal S}^{n\times n})$
such that 
\begin{equation}
u_i>0, \quad
A^{u_i}+\beta_i\in U\qquad \mbox{in}\ \Omega,
\label{A2-1}
\end{equation}
and, for any compact subset $K$ of $\Omega$,
\begin{equation}
u_i\to u, \ \  \beta_i\to 0,
\quad  \mbox{in}\ C^0(K),
\label{A2-2}
\end{equation}
\begin{equation}
F\left(A^{u_i}+\beta_i\right)\to 0\quad  \mbox{in}\ C^0(K).
\label{A2-3}
\end{equation}
\label{DefA1}
\end{Def}

In $\Bbb R^n$, $n\ge 2$, we use ${\cal A}$ to denote the
set of functions $u$ with the following properties:

\noindent {\bf (A1)}\ $u\in C^1(\Bbb R^n), \ u>0$ in $\Bbb R^n$, 
and $\Delta u\le 0$ in $\Bbb R^n$ in the distribution sense.

\noindent{\bf (A2)}\ There exists some $\bar\delta>0$ such that for all
$0<\delta<\bar \delta$, all $x\in \Bbb R^n$, all $\lambda>0$, and
all bounded open set $\Omega$ of $\{y\in \Bbb R^n\
|\ |y-x|>\lambda\}$, 
$$
(1+\delta)u\ge u_{x,\lambda} \ \mbox{in}\ \Omega
\ \mbox{and}\ (1+\delta)u>u_{x,\lambda}\
\mbox{on}\ \partial\Omega\ \mbox{imply}\
(1+\delta)u>u_{x,\lambda}\ \mbox{in}\ \Omega.
$$

\begin{thm} For any $u\in {\cal A}$, there exist
$\bar x\in \Bbb R^n$ and constants $a>0$ and $b\ge 0$ such that
\begin{equation}
u(x)\equiv \left(\frac a{ 1+b|x-\bar x|^2 }\right)^{\frac {n-2}2 },
\qquad \forall\ x\in \Bbb R^n.
\label{B15-1}
\end{equation}
\label{thmB15}
\end{thm}

\begin{rem} If $b=0$, then $u\equiv u(0)$.
\end{rem}

\begin{rem} Theorem \ref{thmB15} can be viewed as an extension 
of the classical result which asserts 
 that positive harmonic functions
in $\Bbb R^n$ are constants.
\end{rem}

\begin{thm} Let $U\subset {\cal S}^{n\times n}$ be an open set
and $F\in C^1(U)\cap C^0(\overline U)$ satisfy 
(\ref{21-1}),
(\ref{C3-1}), 
(\ref{A1-0}), 
(\ref{A3-0}),
(\ref{21-3}), 
(\ref{21-4}), 
and (\ref{A1-1}).
Assume that a positive function $u\in C^1(\Bbb R^n)$ is
a weak solution of
$$
F(A^u)=0\qquad \mbox{in}\ \Bbb R^n,
$$
and satisfies
\begin{equation}
\Delta u\le 0\ \ \mbox{in}\
\Bbb R^n\ \mbox{in the distribution sense}.
\label{A3-1}
\end{equation}
Then
$u\equiv u(0)$ in $\Bbb R^n$.
\label{thmA3}
\end{thm}

\begin{rem} Our notion of weak solutions includes those
arising from rescaling  blowing up solutions.
\end{rem}

\begin{thm} Let $U\subset {\cal S}^{n\times n}$ be a
convex  open set
and $F\in C^1(U)\cap C^0(\overline U)$ satisfy
(\ref{21-1}),
(\ref{C3-1}),
(\ref{A1-0}),
(\ref{A3-0}),
(\ref{21-3}),
(\ref{21-4}),
and (\ref{A1-1}).
Assume  that a positive function $u\in C^{1,1}(\Bbb R^n)$
satisfies (\ref{A3-1}) and
\begin{equation}
F\left(A^u\right)=0,\ \mbox{or equivalently}\
A^u\in \partial U,\  \mbox{almost everywhere}\  \mbox{in}\
\Bbb R^n.
\label{ae}
\end{equation}
Then
$u\equiv u(0)$ in $\Bbb R^n$.
\label{c11}
\end{thm}

\begin{rem} If $U\subset \{M\in {\cal S}^{n\times n}\
|\ Trace (M)\ge 0\}$, then a weak solution or $C^{1,1}$ solution
$u$ of $F(A^u)=0$ automatically satisfies $\Delta u\le 0$.
In particular, this is the case for $(F,U)=(F_k, U_k)$
for all $k$ in all dimensions $n$.
\end{rem}

\begin{rem} It  was proved
by  Chang, Gursky and Yang in
\cite{CGY2} that 
positive $C^{1,1}(\Bbb R^4)$ solutions to
$F_2(A^u)=0$ are constants.
  Aobing Li proved in \cite{Lia} that
   positive $C^{1,1}(\Bbb R^3)$  solutions to
  $F_2(A^u)=0$ are constants, and, for all $k$ and $n$,
   positive $C^3(\Bbb R^n)$
   solutions to $F_k(A^u)=0$ are constants.
 Our proof is completely different.
\end{rem}

We give in the following a notion of weak solutions to
more general equations.

\begin{Def}   Let $U\subset {\cal S}^{n\times n}$ be an open
set, $F\in C^1(U)\cap C^0(\overline U)$, and $h$
 be a continuous function
on an open subset $\Omega$ of $\Bbb R^n$.  A positive function
$u\in C^0(\Omega)$ is said to be a weak solution
of
\begin{equation}
F\left( A^u\right) \le h\qquad\mbox{in}\
\Omega
\label{C1-1}
\end{equation}
if there exist $\{u_i\}$ in $C^2(\Omega)$ and
$\beta_i$ in $C^0(\Omega,
{\cal S}^{n\times n})$ such that, for any
compact subset $K$ of $\Omega$,  (\ref{A2-1}), 
(\ref{A2-2}) hold, and
\begin{equation}
[  F(A^{u_i}+\beta_i) -h] ^+\to 0
\qquad  \mbox{in}\ C^0(K),
\label{C2-3}
\end{equation}
where we have used the notation $w^+:= \max\{w, 0\}$.

Similarly we define that $u$ is a weak solution of 
\begin{equation}
F(A^u)\ge h\qquad\mbox{in}\ \Omega
\label{C2-1}
\end{equation}
by changing $
[  F(A^{u_i}+\beta_i) -h] ^+$ to $[h-F(A^{u_i}+\beta_i)]^+$ in
(\ref{C2-3}).  

We say that $u$ is a weak solution of
$$
F(A^u)=h\qquad \mbox{in}\ \Omega
$$
if it is a weak solution of both (\ref{C1-1}) and (\ref{C2-1}). 
\label{DefC1}
\end{Def}


We also establish the following results concerning degenerate equations
on $\Bbb R^n\setminus\{0\}$.

\begin{thm}
 Let $U\subset {\cal S}^{n\times n}$ be an open set
 and $F\in C^1(U)\cap C^0(\overline U)$ satisfy 
(\ref{21-1}),
(\ref{C3-1}), 
(\ref{A1-0}), 
(\ref{A3-0}),
(\ref{21-3}), 
(\ref{21-4}), 
and (\ref{A1-1}).
 Assume that a positive function $u\in C^1(\Bbb R^n\setminus\{0\})$
  is
  a weak solution  of
$$
  F(A^u)=0\qquad \mbox{in}\ \Bbb R^n\setminus\{0\},
$$
and satisfies
\begin{equation}
\Delta u\le 0\ \ \mbox{in}\
\Bbb R^n\setminus\{0\}\ \mbox{in the distribution sense}.
\label{A3-1new}
\end{equation}
Then
\begin{equation}
u_{x,\lambda}(y)
\le u(y),\qquad \forall\
1<\lambda<|x|, |y-x|\ge \lambda, y\ne 0.
\label{JJJ}
\end{equation}
Consequently,
$u$ is radially symmetric about the origin.
\label{thmA3new}
\end{thm}

\begin{thm}
In addition to the hypotheses on $(F,U)$ in Theorem \ref{thmA3new},
we assume that $U$ is convex.  Assume that
a positive $C^{1,1}(\Bbb R^n\setminus \{0\})$
function $u$ satisfies
(\ref{A3-1new}) and
$$
F\left(A^u\right)=0,\ \mbox{or equivalently}\
A^u\in \partial U,\  \mbox{almost everywhere}\  \mbox{in}\
\Bbb R^n\setminus\{0\}.
$$
Then
(\ref{JJJ}) holds and,
consequently, $u$ is radially symmetric about the origin.
\label{thm12}
\end{thm}

Our proofs of  Theorem \ref{thmA3new} and  Theorem \ref{thm12}
make use of a 
result in a companion paper  
\cite{Li1}.

Let  $\varphi$  be a $C^1$ function
near $1$ satisfying 
$
\varphi(1)=1,
$
and let, for a function $v$, and for $x, y\in \Bbb R^n$,
$\lambda$ close to $1$,
$$
\Phi(v, x, \lambda; y):= \varphi(\lambda) v(x+
\lambda y).
$$

\begin{thm} (\cite{Li1})\ Let
 $\Omega\subset \Bbb R^n$ be a bounded open
 set containing the origin $0$.
 We assume that
 $u\in C^0(\Omega\setminus\{0\})$,
 $v$ is $ C^1$ in some open neighborhood of $\overline \Omega$,
 \begin{equation}
 v>0\qquad\mbox{in}\ \overline \Omega,
 \label{2-1newnew}
 \end{equation}
 \begin{equation}
 \Delta u\le 0  \qquad \mbox{in}\ \Omega\setminus\{0\},
 \label{2-3newnew}
 \end{equation}
 \begin{equation}
 u\ge v\ \mbox{in}\ \Omega\setminus\{0\}.
 \label{2-2newnewnew}
 \end{equation}
Assume that $\varphi$ is as above and
	    \begin{equation}
	    \varphi'(1)v(y)
	    +\nabla v(y)\cdot y>0,\qquad \forall\ y\in \overline \Omega,
	     \label{new}
	      \end{equation}
	       and assume that
	         there exists some $\epsilon_4>0$ such that
		   for any $|x|<\epsilon_4$ and $|\lambda-1|<\epsilon_4$,
		     \begin{equation}
		       \inf_{\Omega\setminus\{0\}}[u-
		         \Phi(v, x, \lambda; \cdot)]=0
			   \ \mbox{implies}\
			    \liminf_{ |y|\to 0} [u(y)-\Phi(v, x, \lambda; y)]=0.
			      \label{gg5znew}
			        \end{equation}
				  Then either
				  \begin{equation}
				  \liminf_{|x|\to 0}[u(x)-v(x)]>0,
				  \label{3-2y}
				  \end{equation}
				  or $u=v=v(0)$ near the origin.
				    \label{thm10}
				      \end{thm}

The paper is organized as follows.  In Section 2, we prove
Theorem \ref{thmB15}.
In Section 3, we give some properties of weak solutions
and $C^{1,1}$ solutions.
In particular, we give comparison principles, see 
Propostion \ref{lemA5-2} and Proposition \ref{prop2},
for weak solutions. 
A crucial ingredient in our proof of
 the comparison principles is Lemma \ref{lemv}, 
 ``the first variation'' of the operator $A^u$.
   Theorem \ref{thmA3} 
and Theorem \ref{c11} are proved in Section 3 by first 
showing  that $u$ belongs to ${\cal A}$
and then showing that the $b$ in (\ref{B15-1}) must be zero.
In Section 4 we first establish some further comparison principles,
Proposition \ref{prop4-1} and Proposition \ref{prop4-2},
which allow the presence of isolated singularities.  Then 
we prove  Theorem \ref{thmA3new} and
Theorem \ref{thm12}.
Proofs of results in this section hevily rely on Theorem \ref{thm10}, a result in 
the companion paper \cite{Li1}.

\section{ Proof of Theorem \ref{thmB15} }

\noindent {\bf Proof of Theorem \ref{thmB15}.}\ By (A1) and the maximum 
principle,
$$
\liminf_{|y|\to\infty} |y|^{n-2}u(y)>0.
$$
As in the proof of lemma 2.1 in \cite{LZ} or \cite{LL1},
for any $x$ in $\Bbb R^n$, there exists 
$\lambda_0(x)>0$ such that
$$
u_{x,\lambda}(y)\le u(y),\qquad \forall\
0<\lambda<\lambda_0(x), |y-x|\ge \lambda.
$$
For any $\delta \in (0, \bar\delta)$, we define 
$$
\bar\lambda_\delta(x):=
\sup\{\mu>0\
|\ u_{x,\lambda}(y)\le (1+\delta)u(y),\
\forall\ 0<\lambda<\mu, |y-x|\ge \lambda\}.
$$
If $\bar\lambda_\delta(x)<\infty$ for some $x$, then
\begin{equation}
u_{x, \bar\lambda_\delta(x) }(y)\le (1+\delta)u(y),
\qquad \forall\ |y-x|\ge \bar\lambda_\delta(x).
\label{A13-0}
\end{equation}

\begin{lem}
If $\bar\lambda_\delta(x)<\infty$
 for some $0<\delta<\bar \delta$ and $x\in \Bbb R^n$,
then
$$
\liminf_{ |y|\to\infty}
|y|^{n-2} \left[
(1+\delta)u(y)-u_{x, \bar\lambda_\delta(x)}(y)\right]=0.
$$
\label{lemA13}
\end{lem}

\noindent{\bf Proof.}\ Suppose the contrary, then for some 
$0<\delta<\bar\delta$, $x\in \Bbb R^n$, $\bar \lambda_\delta(x)<\infty$,
and for some $R>1+\bar \lambda_\delta(x)$ and $\epsilon_1>0$, we have
$$
|y|^{n-2}\left[ (1+\delta)u(y)-u_{x, \bar\lambda_\delta}(y)\right]
\ge 2\epsilon_1,\qquad
\forall\ |y-x|\ge R.
$$
As in the proof of (27) in \cite{LZ}, there exists $\epsilon_2>0$ such that
\begin{equation}
 (1+\delta)u(y)-u_{x, \lambda}(y)\ge 
\frac{\epsilon_2}{|y|^{n-2}},\qquad
\forall\ |\lambda-\bar\lambda_\delta(x)|<\epsilon_2, |y-x|\ge R.
\label{A14-1}
\end{equation}
Since
\begin{equation}
 (1+\delta)u(y)-u_{x, \bar\lambda_\delta(x)}(y)=\delta u(y)>0,\qquad
\forall\ |y-x|=\bar\lambda_\delta(x),
\label{A15-0}
\end{equation}
there exists $0<\epsilon_3<\epsilon_2$ such that
$$
 (1+\delta)u(y)>u_{x, \lambda}(y),\qquad
\forall\ |\lambda -\bar\lambda_\delta(x)|<\epsilon_3,\
 |y-x|=\lambda.
$$

Let
$$
\Omega:=\{ y\in \Bbb R^n\ |\ 
\bar\lambda_\delta(x)<|y-x|<R\}.
$$
We know from (\ref{A13-0}), (\ref{A14-1})
and (\ref{A15-0}) that
$$
(1+\delta)u\ge u_{ x, \bar\lambda_\delta(x) }\
\mbox{in}\ \Omega, \quad 
(1+\delta)u> u_{ x, \bar\lambda_\delta(x) }\
\mbox{on}\ \partial \Omega.
$$
Thus, by (A2), 
\begin{equation}
(1+\delta)u> u_{ x, \bar\lambda_\delta(x) }\
\mbox{on}\ \overline \Omega.
\label{B16-1}
\end{equation}
With (\ref{A14-1}) and (\ref{B16-1}), the moving plane procedure 
can go beyond $\bar\lambda_\delta(x)$, violating the definition 
of $\bar\lambda_\delta(x)$.  Lemma \ref{lemA13}
is established.

\vskip 5pt
\hfill $\Box$
\vskip 5pt

\begin{lem} For all
$0<\delta<\bar\delta$ and for all $x\in \Bbb R^n$,
\begin{equation}
\bar\lambda_\delta(x)^{n-2}
u(x)=(1+\delta)\liminf_{|y|\to\infty}
|y|^{n-2}u(y).
\label{A16-1}
\end{equation}
\label{lemA16-1}
\end{lem}

\noindent{\bf Proof.}\ Let
$\bar\lambda_\delta(\bar x)<\infty$ for
some $0<\delta<\bar \delta$ and  some $\bar x\in \Bbb R^n$.
By Lemma \ref{lemA13},
\begin{equation}
(1+\delta)\liminf_{|y|\to\infty}|y|^{n-2}u(y)
=\lim_{|y|\to \infty}
|y|^{n-2} u_{\bar x, \bar\lambda_\delta(\bar x)}(y)
=\bar\lambda_\delta(\bar x)^{n-2} u(\bar x)<\infty.
\label{A17-1}
\end{equation}

For any $x\in \Bbb R^n$, 
$$
(1+\delta)u(y)\ge u_{x,\lambda}(y),
\qquad\forall\ 0<\lambda<\bar\lambda_\delta(x),\
|y-x|\ge \lambda.
$$
Multiplying the above by $|y|^{n-2}$ and
sending $|y|$ to infinity leads to 
\begin{equation}
(1+\delta)\liminf_{|y|\to\infty}
|y|^{n-2} u(y)\ge
\lambda^{n-2} u(x),\qquad
\forall\ 0<\lambda<\bar\lambda_\delta(x).
\label{A18-1}
\end{equation}

We deduce from (\ref{A18-1}) and (\ref{A17-1}) that
$$
\bar\lambda_\delta(x)^{n-2} u(x)\le
\bar\lambda_\delta(\bar x)u(\bar x)<\infty,
\qquad\forall\ x\in \Bbb R^n.
$$
Switching the roles of $x$ and $\bar x$ leads to (\ref{A16-1})
in the case that $\bar\lambda_\delta$ is not identically
equal to infinity.   On the other hand,
if $\bar\lambda_\delta\equiv \infty$ on $\Bbb R^n$, we send $\lambda$ to $\infty$ in
(\ref{A18-1}) to obtain  (\ref{A16-1}).     Lemma \ref{lemA16-1} is
established.

\vskip 5pt
\hfill $\Box$
\vskip 5pt

If $\bar\lambda_\delta\equiv \infty$
on $\Bbb R^n$ for all $0<\delta<\bar \delta$, then
$$(1+\delta)u(y)\ge u_{x,\lambda}(y),
\quad \forall\ 0<\delta<\bar\delta, x\in \Bbb R^n,
|y-x|\ge \lambda>0.
$$
Sending $\delta$ to $0$ in the above yields
$$
u(y)\ge u_{x,\lambda}(y),
\quad \forall\  x\in \Bbb R^n,
|y-x|\ge \lambda>0.
$$
This implies $u\equiv u(0)$, see  for example   lemma
11.2 in \cite{LZ}.  

We only need to consider the case that for some $0<\delta<\bar\delta$, 
$\bar\lambda_\delta$ is not identically equal to infinity.  According to 
Lemma \ref{lemA16-1} and Lemma \ref{lemA13}, 
$\bar \lambda_\delta(x)<\infty$
for all $x\in \Bbb R^n$,
$$
(1+\delta)u(y)\ge u_{x, \bar\lambda_\delta(x)}(y),\qquad
\forall\ x\in \Bbb R^n, |y-x|\ge \bar\lambda_\delta(x),
$$
and
$$
\liminf_{ |y|\to\infty}
|y|^{n-2}\left[
(1+\delta)u(y)- u_{x, \bar\lambda_\delta(x)}(y)\right]=0,\qquad
\forall\ x\in \Bbb R^n.
$$

For $x\in \Bbb R^n$, let
$$
\phi_\delta^{ (x) }(y):= x+
\frac {  \bar\lambda_\delta(x)^2 (y-x) }{  |y-x|^2  },
\quad \psi(y):=\frac y{ |y|^2 }, \quad
w^{ (x) }:=\left( u_{   \phi_\delta^{  (x)  }  }\right)_\psi
=u_{      \phi_\delta^{  (x)  }  \circ \psi  },
$$
where we have used the notation
$u_\psi:=|J_\psi|^{  \frac {n-2}{2n}  }(u\circ \psi)$ with
$J_\psi$ being the Jacobian of $\psi$.

It is not difficult to see that $w^{  (x)  }$ is $C^1$ near
$0$,
$$
w^{  (x)  }(0)=\bar\lambda_\delta(x)^{n-2}u(x)
=(1+\delta) \liminf_{ |y|\to \infty}
|y|^{n-2} u(y)=\liminf_{ |y|\to 0} (1+\delta) u_\psi(y)>0,
$$
$$
w^{  (x)  }\le (1+\delta)u_\psi\qquad
\mbox{in}\ B_{\delta(x)}\setminus\{0\}\
\mbox{for some}\  \delta(x)>0.
$$
Following the proof of theorem 1.3 
in \cite{LL3} (see also \cite{LL2} and \cite{LL4}),  we obtain
$$
\nabla w^{  (x)  }(0)=\nabla w^{  (0)  }(0),
\qquad \forall\ x\in \Bbb R^n
$$
and then
$$
u(y)\equiv \left(    \frac{ \alpha_\delta^{ \frac 2{n-2}  }  }
{  d+|y-\bar x|^2  }   \right)^{  \frac{n-2}2  },
\qquad\forall\ y\in \Bbb R^n,
$$
where
$\bar x\in \Bbb R^n, d>0$ and 
$$
\alpha_\delta= (1+\delta) 
\liminf_{  |y|\to \infty }  |y|^{n-2}u(y)>0.
$$
Theorem \ref{thmB15} is established.

\section{Properties of weak solutions and the proof
 of Theorem \ref{thmA3} and Theorem \ref{c11}}

We start with some
properties
of weak solutions.

\begin{lem} Let  $U\subset {\cal S}^{n\times n}$ be an open
set satisfying (\ref{C3-1}).
  Assume that for a positive
$C^2$ function $u$ in
some open subset $\Omega$ of $\Bbb R^n$,
 there exist $\{u_i\}\subset C^2(\Omega)$
and $\{\beta_i\}
\subset C^0(\Omega,
{\cal S}^{n\times n})$
such that (\ref{A2-1}) and (\ref{A2-2}) hold on any compact subset $K$ of
$\Omega$.
Then
\begin{equation}
A^u\in \overline U\qquad\mbox{in}\ \Omega.
\label{C3-2}
\end{equation}
\label{lemC3}
\end{lem}

\noindent{\bf Proof.}\
For any $\bar x\in \Omega$, we fix some $\delta>0$ such that
$B_{2\delta}(\bar x)\subset \Omega$.
Consider, for small $\epsilon>0$,
$$
u^\epsilon(y):= u(y)-\frac \epsilon 2 |y-\bar x|^2,\qquad
y\in B_\delta(\bar x).
$$
We know from (\ref{A2-2}) that 
$$
u_i(\bar x)=u^\epsilon(\bar x)+\circ(1),
$$
$$
u_i(y)=u^\epsilon(y)+\frac 12
\epsilon\delta^2 +\circ(1),\qquad
y\in \partial B_\delta(\bar x),
$$
$$
u_i(y)\ge u^\epsilon(y)+\circ(1), \qquad y\in  B_\delta(\bar x),
$$
where $\circ(1)\to 0$ as $i\to \infty$, uniform in
$y$ and $\epsilon$.

It is easy to see that for some $a_i^\epsilon=1+\circ(1)$,
$y_i^\epsilon\in B_\delta(\bar x)$,
\begin{equation}
u_i\ge a_i^\epsilon u^\epsilon\qquad
\mbox{on}\ B_\delta(\bar x),
\label{C5-1}
\end{equation}
\begin{equation}
u_i(y_i^\epsilon)=a_i^\epsilon u^\epsilon(y_i^\epsilon).
\label{C5-2}
\end{equation}
Passing to a subsequence in (\ref{C5-2}),
$y_i^\epsilon\to \bar y^\epsilon$,
$
u(\bar y^\epsilon)=u(\bar y^\epsilon)-\frac
  \epsilon 2 |\bar y^\epsilon-\bar x|^2$.  So
$$
\lim_{i\to \infty} y_i^\epsilon=\bar x.
$$
By (\ref{C5-1}) and (\ref{C5-2}),
\begin{equation}
A^{u_i}(y_i^\epsilon)\le A^{  a_i^\epsilon u^\epsilon}(y_i^\epsilon).
\label{C6-1}
\end{equation}
Thus, by (\ref{A2-1}) and (\ref{C3-1}),
$ A^{  a_i^\epsilon u^\epsilon}(y_i^\epsilon) +\beta_i(y_i^\epsilon)\in
\overline U$.  Sending $i$ to $\infty$, we have, using (\ref{A2-2}), 
$
A^{ u^\epsilon}(\bar x)\in \overline U.
$
Sending $\epsilon$ to $0$, we have
$A^u(\bar x)\in \overline U.
$
Lemma \ref{lemC3} is established.

\vskip 5pt
\hfill $\Box$
\vskip 5pt

\begin{lem}  Let  $U\subset {\cal S}^{n\times n}$ be an open
set satisfy (\ref{C3-1}) and  $F\in C^1(U)\cap C^0(\overline U)$
satisfy (\ref{21-4}).   Assume that
for a  positive
$C^2$ function $u$ in
some open set  $\Omega$ of $\Bbb R^n$,
there exist $\{u_i\}\subset C^2(\Omega)$ and
$\{\beta_i\}\subset C^0(\Omega, {\cal S}^{n\times n})$ such that
(\ref{A2-1}) and (\ref{A2-2}) hold for any
compact subset $K$ of $\Omega$, and, for some
$h\in C^0(\Omega)$,
\begin{equation}
[h-F(A^{u_i}+\beta_i)]^+\to 0,\qquad \mbox{in}\ C^0_{loc}(\Omega).
\label{KK}
\end{equation} 
  Then (\ref{C3-2}) holds and $u$ is 
a classical solution of (\ref{C2-1}).
\label{lemC7-1}
\end{lem}

\noindent {\bf Proof.}\  We know from Lemma \ref{lemC3} that 
(\ref{C3-2}) holds.  Following the proof
of Lemma \ref{lemC3} from the beginning until (\ref{C6-1}).
Then, by (\ref{A2-1}),  
 (\ref{C3-1}), (\ref{C6-1}) and (\ref{21-4}),
$$
A^{a_i^\epsilon u^\epsilon }(y_i^\epsilon)
+\beta_i(y_i^\epsilon)\in \overline U,
$$
and
\begin{equation}
F\left(  A^{u_i}(y_i^\epsilon)+\beta_i(y_i^\epsilon)\right)\le
F\left(  A^{a_i^\epsilon u^\epsilon }(y_i^\epsilon)
+\beta_i(y_i^\epsilon)\right).
\label{9A}
\end{equation}
Since
$$
\lim_{i\to \infty}
\left(  A^{a_i^\epsilon u^\epsilon }(y_i^\epsilon)
+\beta_i(y_i^\epsilon)\right)=
A^{u^\epsilon}(\bar x),
$$
we have $A^{u^\epsilon}(\bar x)\in \overline U$
and, using the continuity of $F$ on $\overline U$,
$$
\lim_{i\to \infty}
F\left(  A^{a_i^\epsilon u^\epsilon }(y_i^\epsilon)
+\beta_i(y_i^\epsilon)\right)=
F\left(
A^{u^\epsilon}(\bar x)
\right).
$$
Sending $i$ to $\infty$ in (\ref{9A}) leads to, in view of
(\ref{KK}),
$$
F\left(  A^{u^\epsilon }(\bar x)\right) \ge h(\bar x).
$$
Sending $\epsilon$ to $0$, we obtain 
$$
F\left(  A^u(\bar x) \right)\ge h(\bar x).
$$
Lemma \ref{lemC7-1} is established.

  \vskip 5pt
    \hfill $\Box$
      \vskip 5pt

\begin{lem}
 Let $U\subset {\cal S}^{n\times n}$ be an open
 set satisfy (\ref{C3-1}), $F\in C^1(U)\cap C^0(\overline U)$
 satisfy (\ref{21-4}).  Assume that
 for a positive $C^2$ function $u$ in
 some open set  $\Omega$ of $\Bbb R^n$,
there exist $\{u_i\}\subset C^2(\Omega)$ and
$\{\beta_i\}\subset C^0(\Omega, {\cal S}^{n\times n})$
such that (\ref{A2-1}),  (\ref{A2-2})
and
\begin{equation}
\sup_i \sup_K |\nabla^2u_i|<\infty
\label{A2-0}
\end{equation}
hold  for
any compact subset $K$ of $\Omega$, and, for some
$h\in C^0(\Omega)$,
\begin{equation}
\left[F(A^{u_i}+\beta_i)-h\right]^+
\to 0, \qquad \mbox{in}\
C^0_{loc}(\Omega).
\label{KKK}
\end{equation}
    Then (\ref{C3-2}) holds and $u$ is
 a classical solution of (\ref{C1-1}).
 \label{lemC9}
 \end{lem}

\noindent {\bf Proof.}\ We know from Lemma \ref{lemC3} that (\ref{C3-2})
holds.  
For any $\bar x\in \Omega$, we fix some 
$\delta>0$ such that
$B_{2\delta}(\bar x)\subset \Omega$.
Consider,
for small $\epsilon>0$,
$$
u^\epsilon(y):= u(y)+\frac \epsilon 2 |y-\bar x|^2,
\qquad y\in B_\delta(\bar x).
$$
Arguing as in the proof of Lemma \ref{lemC3}, we find
$a_i^\epsilon=1+\circ(1)$, $y_i^\epsilon\to \bar x$,
\begin{equation}
A^{u_i}(y_i^\epsilon)\ge A^{  a_i^\epsilon u^\epsilon }(y_i^\epsilon).
\label{C11-1}
\end{equation}
Clearly, there exist $\alpha_\epsilon>0$
and $\gamma_i>0$ satisfying
$\alpha_\epsilon\to 0$ as $\epsilon\to 0$ and $\gamma_i\to 0$ as $i\to \infty$,
such that
\begin{equation}
A^{  a_i^\epsilon u^\epsilon }(y_i^\epsilon)
+\beta_i(y_i^\epsilon)\ge
A^u(\bar x)-\alpha_\epsilon I-\gamma_i I.
\label{C11-2}
\end{equation}
We already know that $A^u(\bar x)\in 
\overline U$.  So,
by (\ref{C11-1}), (\ref{C11-2}), (\ref{C3-1}) and
(\ref{21-4}),
$$
  A^{u_i}(y_i^\epsilon)
+\beta_i(y_i^\epsilon)
+\alpha_\epsilon+\gamma_i I
\in \overline U
$$
 and
$$
F\left( A^u(\bar x)\right)\le
F\left(  A^{u_i}(y_i^\epsilon)+\beta_i(y_i^\epsilon)
+\alpha_\epsilon I
+\gamma_i I\right).
$$
Because of (\ref{A2-0}), (\ref{A2-2}) and the positivity and the continuity
 of
$u$, 
$  A^{u_i}(y_i^\epsilon)+\beta_i(y_i^\epsilon)
+\alpha_\epsilon I
+\gamma_i I$ remain bounded.  Thus, by the continuity of $F$,
$$
F\left( A^u(\bar x)\right)\le
F\left(  A^{u_i}(y_i^\epsilon)+\beta_i(y_i^\epsilon)\right)
+\circ(1)+\circ_\epsilon(1),
$$
where $\circ_\epsilon(1)\to 0$ as $\epsilon\to 0$, uniform in $i$, and
$\circ(1)\to 0$ as $i\to \infty$, uniform in $\epsilon$.
Sending $i$ to $\infty$ and then $\epsilon$ to $0$, we obtain,
using (\ref{KKK}), 
\begin{equation}
F\left( A^u(\bar x)\right)\le
h(\bar x).
\label{HH1}
\end{equation}
Lemma \ref{lemC9} is established.

  \vskip 5pt
    \hfill $\Box$
      \vskip 5pt

\begin{lem}
In Lemma \ref{lemC9}, we drop assumption (\ref{A2-0}) but
add $A^u\in U$ in $\Omega$.
Then  $u$ is
 a classical solution of (\ref{C1-1}).
  \label{lemC10}
   \end{lem}

\noindent {\bf Proof.}\
Follow the proof of Lemma \ref{lemC9} until 
(\ref{C11-2}).  Since we know that $A^u(\bar x)\in U$,
the right hand side of (\ref{C11-2})
 is also in $U$ for large $i$ and small $\epsilon$.
Thus
$$
F\left(A^u(\bar x)-\alpha_\epsilon I-\gamma_i I\right)
\le F\left( A^{u_i}(y_i^\epsilon)
+\beta_i(y_i^\epsilon) \right).
$$
Sending $i$ to $\infty$ and then
$\epsilon$ to $0$, we obtain 
(\ref{HH1}). Lemma \ref{lemC10} is established.

  \vskip 5pt
      \hfill $\Box$
            \vskip 5pt

\begin{lem}  Let $\Omega\subset \Bbb R^n$ be
an open set. Assume that a positive function 
$u\in C^0(\Omega)$ is a weak solution of (\ref{A1-2}).
  Then, for any
constant $b>0$ and for any $x\in \Bbb R^n$,
the function $
\displaystyle{
v(y):=b^{ \frac{n-2}2}u(x+by)
}$
 is  a weak solution of
$$
F(A^v)=0\qquad \mbox{in}\ \widehat \Omega:=\{y\in \Bbb R^n\
|\ x+by\in \Omega\}.
$$
\label{lemA5-1}
\end{lem}

\noindent{\bf Proof.}\ It is obvious.

\begin{lem}
Let $\Omega\subset \Bbb R^n$ be
an open set. Assume that  a positive function
$u\in C^0(\Omega)$ is a weak solution of (\ref{A1-2}).
  Then,
for any $x\in \Bbb R^n$ and $\lambda>0$,
$u_{x,\lambda}$ is a weak solution of
$$
F\left( A^{ u_{x,\lambda} }\right)=0\qquad
\mbox{in}\ \Omega_{x,\lambda}
:=
\left\{ y\in \Bbb R^n\ \bigg|\
x+ \frac{ \lambda^2(y-x)  }{  |y-x|^2 }
\in \Omega\right\}.
$$
\label{lemA10-1}
\end{lem}

\noindent{\bf Proof.}\ This follows from the conformal invariance
of the operator $F(A^u)$, see for example line 9 on page 
1431 of \cite{LL1}.

\vskip 5pt
\hfill $\Box$
\vskip 5pt

The following is a comparison principle for weak solutions.
\begin{prop} Let $U\subset {\cal S}^{n\times n}$ be an open set satisfying
(\ref{C3-1}) and (\ref{A1-0}), $F\in C^1(U)\cap C^0(\overline U)$ satisfy
(\ref{21-4}) and (\ref{A1-1}),  
and let  $\Omega\subset \Bbb R^n$ 
be a bounded open set,  
 $u, v\in C^0(\overline \Omega)$ satisfy
\begin{equation}
u, v>0,\qquad \mbox{in}\ \overline \Omega,
\label{H1}
\end{equation}
and
\begin{equation}
u>v\qquad\mbox{on}\ \partial \Omega.
\label{A6-2}
\end{equation}
We assume that there exist
$\{\beta_i\}, \{\tilde \beta_i\}\subset C^0(\Omega, {\cal S}^{n\times n})$
and  positive functions $\{u_i\}, \{v_i\}\subset 
C^2(\Omega)$ such that, for any
compact subset $K$ of $\Omega$,
\begin{equation}
A^{u_i}+\beta_i\in U, \ A^{v_i}
+\tilde \beta_i\in U,\qquad \mbox{in}\ \Omega,
\label{H2}
\end{equation}
\begin{equation}
u_i\to u,\ \ v_i\to v, \ \beta_i\to 0,\
\tilde \beta_i \to 0,  \qquad
\mbox{in}\ C^0(K),
\label{H3}
\end{equation}
and
\begin{equation}
F\left(A^{v_i}+\tilde \beta_i\right)\to 0\quad  \mbox{in}\ C^0(K).
\label{H5new}
\end{equation}
Then
\begin{equation}
u>v \qquad\mbox{on}\ \overline \Omega.
\label{A6-3}
\end{equation}
\label{lemA5-2}
\end{prop}

To prove Proposition \ref{lemA5-2},
we need to produce appropriate approximations to the $\{u_i\}$.
This is achieved by studying ``the first variation''
of the operator $A^u$.

Writing 
$$
w=u^{-\frac 2{n-2} },
$$
we have
$$
A^u=A_w:= w\nabla ^2w -\frac 12 |\nabla w|^2 I.
$$

\begin{lem} Let $\Omega$ be a bounded open set in $\Bbb R^n$,
$w\in C^2(\Omega)$ satisfy, for some constant $c_1>0$,
$$
w\ge c_1,\qquad \mbox{in}\ \Omega,
$$
and
let
\begin{equation}
\varphi(y)=e^{ \delta |y|^2}.
\label{w0}
\end{equation}
Then there  exists some constant $\delta>0$,
depending only on $\sup\{|y|\ |\ y\in \Omega\}$,
and there exists
$\bar\epsilon$, depending only on $\delta$, $c_1$ and
 $\sup\{|y|\ |\ y\in \Omega\}$,
such that
for any  $0<\epsilon<\bar \epsilon$,
\begin{equation}
A_{ w+\epsilon\varphi}\ge
\left(1+\epsilon \frac \varphi w\right)A_w
+\frac {\epsilon \delta}2 \varphi w I\qquad
\mbox{in}\ \Omega.
\label{w6}
\end{equation}
\label{lemv}
\end{lem}

\noindent {\bf Proof of Lemma \ref{lemv}.}\
Let $\varphi$ be  a fixed function, a computation gives
$$
A_{w+\epsilon\varphi}=
A_w+\epsilon\left\{
w\nabla^2\varphi+\varphi \nabla^2w-\nabla w\cdot \nabla \varphi
I\right\}
+\epsilon^2A_\varphi.
$$
Replacing $\nabla^2 w $ by $\displaystyle{
w^{-1} \left( A_w+\frac 12 |\nabla w|^2I\right)}$
in the above, we have
\begin{equation}
A_{w+\epsilon\varphi}=
(1+\epsilon\frac \varphi w)A_w
+\epsilon\left\{ w \nabla^2 \varphi+
\frac {  |\nabla w|^2 }{ 2w} \varphi I
-\nabla w\cdot \nabla \varphi I\right\}
+\epsilon^2A_\varphi.
\label{w3}
\end{equation}

For the $\varphi$ in 
(\ref{w0}),
$$
\nabla \varphi(y)=2\delta \varphi(y)y,
\quad \nabla ^2 \varphi(y)
= 2\delta \varphi(y)I+4\delta^2 \varphi(y) y\otimes y.
$$
It follows that 
\begin{eqnarray}
&& w \nabla^2 \varphi+
 \frac {  |\nabla w|^2 }{ 2w} \varphi I
 -\nabla w\cdot \nabla \varphi I
 \nonumber\\
 &\ge& \varphi\left\{ 2\delta w+  \frac {  |\nabla w|^2 }{ 2w} 
 -2\delta \nabla w\cdot y\right\}I
 \ge   \varphi\left\{ 
 2\delta w+  \frac {  |\nabla w|^2 }{ 2w}
 - \left( \frac{ |\nabla w|}{ \sqrt{4w} }\right)
  \left( 2\delta\sqrt{4w}|y|\right)\right\}I \nonumber\\
  &\ge &   \varphi\left\{
   2\delta w+  \frac {  |\nabla w|^2 }{ 4w}
   -4\delta^2 w|y|^2\right\}I. \nonumber
 \label{w4}
 \end{eqnarray}
It is clear that there exists 
 $\delta>0$, depending only
on $\sup\{ |y|\ |\ y\in\Omega\}$, such that
$$
w \nabla^2 \varphi+
 \frac {  |\nabla w|^2 }{ 2w} \varphi I
  -\nabla w\cdot \nabla \varphi I
  \ge \delta \varphi wI+  \frac {  |\nabla w|^2 }{ 4w}\varphi
  I.
$$
For this $\delta$, there exists $\bar\epsilon>0$, depending only on
$\delta $,  $c_1$ and $\sup\{ |y|\ |\ y\in\Omega\}$,
 such that
for all $0<\epsilon<\bar\epsilon$,
$$
A_{ w+\epsilon\varphi}\ge \left(1+\epsilon \frac \varphi w\right)A_w
+ \epsilon \delta  \varphi w I+\epsilon^2 A_\varphi\ge
\left(1+\epsilon \frac \varphi w\right)A_w
+\frac {\epsilon \delta}2 \varphi w I,
$$
i.e.
$$
A^{ (u^{- \frac 2{n-2}}
+\epsilon \varphi )^{ -\frac {n-2}2 }  }
\ge
\left(1+\epsilon \frac \varphi w\right)A_w
+\frac {\epsilon \delta}2 \varphi w I.
$$
Lemma \ref{lemv} is established.

\vskip 5pt
\hfill $\Box$
\vskip 5pt

\noindent{\bf Proof of Proposition \ref{lemA5-2}.}\
Since shrinking $\Omega$ slightly will not affect
(\ref{A6-2}), we may assume without loss of generality that 
(\ref{H3}) and (\ref{H5new}) hold
with $K$ replaced by $\overline \Omega$ --- from now on these
equations will be understood in this sense.

We prove (\ref{A6-3}) by contradiction argument. 
 Suppose the contrary,
 then there exists some $\bar x\in \Omega$ such that
$$
u(\bar x)\le v(\bar x).
$$
It is clear that there exist $0<a\le 1$ and $\bar y\in \Omega$ such that
\begin{equation}
u\ge av\qquad \mbox{in}\ \Omega,
\label{A7-1}
\end{equation}
\begin{equation}
u(\bar y)=av(\bar y).
\label{A7-2}
\end{equation}

 Let $\delta$ and
 $\varphi$ be as
 in Lemma \ref{lemv},
 $$
 \epsilon_i:= \sup _{y\in \Omega} \sqrt { \|\beta_i(y)\| }
 \to 0,
 $$
 and
 $$
 \widehat u_i:=
 \left( u_i^{ -\frac 2{n-2} }+\epsilon_i\varphi\right)^{ -\frac {n-2}2 }.
 $$
 By (\ref{w6}),
\begin{equation} 
 A^{ \widehat u_i }
 \ge \left( 1+\epsilon_i \frac \varphi{w_i}\right)A^{u_i}
 +\frac {\epsilon_i\delta}2 \varphi  w_i I,
\label{large}
\end{equation}
 where $ w_i= u_i^{ -\frac 2{n-2} }$.

It is easy to see, using (\ref{A7-1}), (\ref{A7-2}), (\ref{A6-2}),
 and the convergence of  $\widehat
u_i$ to $u$ and $v_i$ to $v$, that
for some $\bar \delta>0$ and some large integer
$\bar I$, there exist,  for $i, j\ge \bar I$,
$a_{ij}\in (\frac a2, \frac {3a}2)$ and $y_{ij}\in 
\{y\in
\Omega\ |\ dist(y, \partial \Omega)>\bar \delta\}$ such that
$$
\widehat u_i\ge a_{ij} v_j\qquad
\mbox{in}\ \Omega, 
$$
$$
\widehat u_i(y_{ij})=a_{ij}v_j(y_{ij}).
$$
It follows that
$$
\nabla \left( a_{ij}^{-1}
\widehat u_i\right)(y_{ij})=\nabla v_j(y_{ij}),
\qquad \left(  \nabla^2 (a_{ij}^{-1} \widehat u_i)(y_{ij})\right)
\ge  \left(  (\nabla^2v_j(y_{ij})\right),
$$
and
$$
A^{a_{ij}^{-1}\widehat u_i}(y_{ij})\le A^{  v_j}(y_{ij}).
$$
Thus, in view of (\ref{large}) and the definition of
$\epsilon_i$,
for some $\tilde I\ge \bar I$, and for all $i, j\ge \tilde I$,
\begin{eqnarray*}
A^{  v_j}(y_{ij})+\tilde \beta_j(y_{ij})
&\ge &
(a_{ij})^{ \frac 4{n-2} }A^{u_i}(y_{ij})
+\tilde \beta_j(y_{ij})
\\
&\ge& (a_{ij})^{ \frac 4{n-2} }
 \left( 1+\epsilon_i \frac \varphi{w_i}\right)A^{u_i}
 + (a_{ij})^{ \frac 4{n-2} } \frac {\epsilon_i\delta}2 \varphi  w_i I
 +\tilde \beta_j(y_{ij})\\
 &\ge&
  (a_{ij})^{ \frac 4{n-2} }
   \left( 1+\epsilon_i \frac \varphi{w_i}\right)
   \left( A^{u_i}+\beta_i \right)(y_{ij})
    + (a_{ij})^{ \frac 4{n-2} } \frac {\epsilon_i\delta}4 \varphi  w_i I
     +\tilde \beta_j(y_{ij}).
 \end{eqnarray*}
Fixing $i=\tilde I$, we have,
for large $j$,
$$
A^{  v_j}(y_{\tilde Ij})+\tilde \beta_j(y_{\tilde Ij})
\ge 
  (a_{\tilde Ij})^{ \frac 4{n-2} }
     \left( 1+\epsilon_{\tilde I} \frac \varphi{w_{\tilde I}}\right)
       \left( A^{u_{\tilde I}}+\beta_{\tilde I} \right)(y_{\tilde Ij})
          + (a_{\tilde Ij})^{ \frac 4{n-2} } \frac {\epsilon_{\tilde I}
	  \delta}8 \varphi  w_{\tilde I} I.
	  $$

By (\ref{H2}),  (\ref{C3-1}), (\ref{A1-0}),  (\ref{21-4}) and
(\ref{A1-1}),
\begin{eqnarray*}
&&
F\left(   A^{  v_j}(y_{\tilde Ij})+\tilde \beta_j(y_{\tilde Ij})
\right)\\
&\ge&
F\left(   (a_{\tilde Ij})^{ \frac 4{n-2} }
     \left( 1+\epsilon_{\tilde I} \frac \varphi{w_{\tilde I}}\right)
            \left( A^{u_{\tilde I}}+\beta_{\tilde I} \right)(y_{\tilde Ij})
	    \right)
	   \\
	   &\ge&
\min_{y\in \Omega,
dist(y, \partial \Omega)\ge \bar \delta}
\min_{ (\frac a2)^{ \frac 4{n-2} }
	    \le b\le (\frac{3a}2)^{ \frac 4{n-2} }
	    }
	    F\left(b \left( 1+\epsilon_{\tilde I}
	    \frac \varphi{w_{\tilde I}}\right)
       \left( A^{u_{\tilde I}}+\beta_{\tilde I} \right)(y)
			            \right)>0
\end{eqnarray*}
 Sending $j$ to $\infty$ leads to,
 in view of 
 (\ref{H5new}),  that
 $$
0=\lim_{j\to \infty} 
 F\left(   A^{  v_j}(y_{\tilde Ij})+\tilde \beta_j(y_{\tilde Ij})
 \right)
 >0.
$$ 
  Impossible.  Proposition \ref{lemA5-2} 
   is established.

  \vskip 5pt
  \hfill $\Box$
  \vskip 5pt

\begin{rem} If we further assume in Proposition \ref{lemA5-2} that
$$
\sup_i \sup_K  \left(
 |\nabla^2 u_i|+   |\nabla^2v_i|\right)
 <\infty,
 $$
 then modification of the proof
of Proposition \ref{prop2} yields a somewhat different proof.
One observation is needed: In addition to (\ref{vv4}), we have
$D^2u_{i, \epsilon}\ge D^2 v_{i, \epsilon}$ on
$S_{i, \epsilon}$.  So, for small $s$, using (\ref{C3-1}) and the
continuity of $D^2u_{i,\epsilon}$ and 
$D^2 v_{i,\epsilon}$, we still have
(\ref{vv3}).
\end{rem}

\noindent{\bf Proof of Theorem \ref{thmA3}.}\
We first prove that $u\in {\cal A}$.  We only need to verify property (A2)
 since property (A1)
has already been assumed. Let $\Omega$ be as in the statement of (A2),
then, by
 Lemma \ref{lemA5-1} and Lemma \ref{lemA10-1},
both $u_{x,\lambda}$ 
and $(1+\delta)u$ are  weak solutions of (\ref{A1-2}).
Thus, by Proposition \ref{lemA5-2}, 
(A2) is satisfied.  So we have proved that $u\in {\cal A}$.  
By Theorem \ref{thmB15}, (\ref{B15-1}) holds for some $a>0, b\ge 0$ and
$\bar x\in \Bbb R^n$.  We only need to prove that $b=0$.  Suppose that
$b>0$, then a computation
gives, for some positive constant $\beta$,
\begin{equation}
A^u\equiv \beta I\qquad\mbox{in}\ \Bbb R^n.
\label{B6-1}
\end{equation}

Since $u$ is a weak solution of (\ref{A1-2}), 
the hypotheses of Lemma \ref{lemC3} with
$h\equiv 0$
are satisfied, and therefore, according to Lemma \ref{lemC3},
\begin{equation}
A^u\in \overline U\qquad \mbox{in}\ \Bbb R^n.
\label{III}
\end{equation}
By (\ref{B6-1}), (\ref{III}) and (\ref{A3-0}),
$$
A^u\equiv \beta I\in U\qquad \mbox{in}\ \Bbb R^n.
$$
Thus, by (\ref{A1-1}),
$$
F(A^u)>0\qquad\mbox{in}\ \Bbb R^n,
$$
On the other hand, since $A^u\in U$, we apply 
 Lemma \ref{lemC10} to obtain
$$
F(A^u)\le 0.
$$
Impossible. 
  We have 
proved that $b=0$ in (\ref{B15-1}) and
therefore $u\equiv u(0)$.
Theorem \ref{thmA3} is established.

 \vskip 5pt
   \hfill $\Box$
     \vskip 5pt

To prove  Theorem \ref{c11} we need the following 
comparison principle for $C^{1,1}$ solutions.
\begin{prop} Let $U\subset {\cal S}^{n\times n}$ be a
convex  open set satisfying
(\ref{C3-1}) and (\ref{A1-0}), $F\in C^1(U)\cap C^0(\overline U)$ satisfy
(\ref{21-4}) and (\ref{A1-1}),
and let  $\Omega\subset \Bbb R^n$
be a bounded open set,
 $u, v\in C^{1,1}_{loc}(\Omega)$ satisfy
 (\ref{H1}), (\ref{A6-2}) and
 \begin{equation}
 F\left(A^u\right)\ge 0, \ \ \  F\left(A^v\right)=0,\qquad
 \mbox{almost everywhere in}\ \Omega.
 \label{convex1}
 \end{equation}
 Then (\ref{A6-3}) holds.
 \label{prop2}
 \end{prop}

\noindent{\bf Proof.} 
We prove it by contradition argument. We assume that $
\displaystyle{  \min_\Omega (u-v)\le 0}$.
Let $\varphi$ be as in (\ref{w0}) for some fixed  small $\delta>0$, and
let
$$
w:=u^{-\frac 2{n-2}}, \
w_\epsilon:=w+\epsilon\varphi,\
u_\epsilon:= (w_\epsilon)^{-\frac {n-2}2}.
$$
Using Lemma \ref{lemv} and (\ref{A6-2}), we can find some fixed small 
positive constants $\epsilon$ and $\epsilon_1$ such that
\begin{equation}
A^{u_\epsilon}\ge \left( 1+ \epsilon\frac \varphi{ w}\right)
A^u+ 3\epsilon_1 I,
\label{vv1}
\end{equation}
and
$$
u_\epsilon> v\qquad\mbox{on}\ \partial \Omega.
$$
Since $A^u\in \partial U$ a.e. in $\Omega$,
we have, using (\ref{C3-1}) and the openness of $U$,
\begin{equation}
A^{u_\epsilon}+M\in U\ 
\mbox{a.e. in}\ \Omega,\qquad \forall\ M\in {\cal S}^{n\times n},
\ \|M\|<2\epsilon_1.
\label{vv2}
\end{equation}
By the contradiction hypothesis, $u\le v$ somewhere
in $\Omega$, so
there exists $a_\epsilon\in (0, 1]$
such that
$$
u_\epsilon\ge v_\epsilon:= a_\epsilon v\qquad\mbox{in}\ \Omega,
$$
$$
u_\epsilon>  v_\epsilon\qquad\mbox{on}\ \partial \Omega.,
$$
and
$$
S_\epsilon:=\{x\in \Omega\ |\ u_\epsilon(x)=v_\epsilon(x)\}\ne \emptyset.
$$
Clealy, 
\begin{equation}
u_\epsilon=v_\epsilon, \ \nabla u_\epsilon=\nabla v_\epsilon,\qquad
\mbox{on}\ S_\epsilon.
\label{vv4}
\end{equation}
Recall that $\epsilon$ has been fixed.
Let
$$
O_s:=\{x\in \Omega\ |\
dist(x, S_\epsilon)<s\}.
$$
By (\ref{vv4}),
$$
A^{tu_\epsilon+(1-t)v_\epsilon}=tA^{u_\epsilon}+(1-t)A^{v_\epsilon}\qquad
\mbox{on}\ S_\epsilon.
$$
So, 
in view  
 the convexity and the openness of 
$U$, there exists a function $\bar t(s)$,
$\bar t(s)\to 0^+$ as $s\to 0^+$, such that
\begin{equation}
A^{tu_\epsilon+(1-t)v_\epsilon}\in U
\ \mbox{a.e. in}\ O_s,\qquad
\forall\ \bar t(s)<t\le 1.
\label{vv3}
\end{equation}
Note that we have used, in deriving (\ref{vv3}), 
the fact that $A^u$ is linear in $\nabla^2 u$
and both $u$ and $\nabla u$ are continuous.

By (\ref{vv2}), there exists some $\epsilon_2>0$, independent of
$s$, such that
$$
dist\left(A^{tu_\epsilon+(1-t)v_\epsilon}, \partial U\right)>\epsilon_2,\qquad
\forall\ 1-\epsilon_2\le t\le 1.
$$
Thus, by (\ref{21-4}), there exists some
$\epsilon_3>0$,  independent of
$s$, such that
\begin{equation}
F\left(A^{tu_\epsilon+(1-t)v_\epsilon}\right)\ge \epsilon_3,
\ \ \left(F_{ij}\left(A^{tu_\epsilon+(1-t)v_\epsilon}\right)\right)
\ge \epsilon_3I,\qquad
\forall\ 1-\epsilon_2\le t\le 1.
\label{vv5}
\end{equation}
Using the mean value theorem, in view of (\ref{vv5}), we have
\begin{eqnarray}
&&\epsilon_3-F\left(A^{\bar t(s)u_\epsilon+(1-\bar t(s))v_\epsilon}\right)
\nonumber\\
&\le & F\left(A^{u_\epsilon}\right)-
F\left(A^{\bar t(s)u_\epsilon+(1-\bar t(s)v_\epsilon}\right)
= \int_{ \bar t(s) }^1
\left\{  \frac {d}{dt} F\left(A^{tu_\epsilon+(1-t)v_\epsilon}\right) \right\}
dt\nonumber\\
&=:&  -\left(\int_{ \bar t(s) }^1 a_{ij}(x, t)dt\right)
\partial_{ij}(u_\epsilon-v_\epsilon)+b_i(x)\partial_i(u_\epsilon-v_\epsilon)
+c(x)(u_\epsilon-v_\epsilon), \label{vv6} 
\end{eqnarray}
where $a_{ij}(\cdot, t), b_i, c$ are bounded in $L^\infty$ norm, and,
in view of  (\ref{21-4}) and (\ref{vv5}),
$$
 \left(\int_{ \bar t(s) }^1 a_{ij}(x, t)dt\right)
 \ge \epsilon_4I
 $$
 for some $\epsilon_4>0$ independent of $s$.
In view of  (\ref{vv6}), we can find
some small $\bar s>0$ such that
$$
0<\frac 12 \epsilon_3\le 
-a_{ij}\partial_{ij}(u_\epsilon-v_\epsilon)+b_i(x)\partial_i(u_\epsilon-v_\epsilon)
+c(x)(u_\epsilon-v_\epsilon),\qquad \mbox{a.e. in }\ O_{\bar s},
$$
where $a_{ij}, b_i, c$ are in $L^\infty(O_{\bar s})$ and
$\displaystyle{
\left(a_{ij}\right)\ge \epsilon_4I}$ a.e. in $O_{\bar s}
$.
We know that
$$
u_\epsilon-v_\epsilon=0\ \mbox{on}\
S_\epsilon\subset O_s,\ 
u_\epsilon-v_\epsilon\ge 0\ \mbox{in}\ O_s,
u_\epsilon-v_\epsilon>0\ \mbox{near}\ \partial O_{\bar s}.
$$
But this violates the local maximum principle,
see theorem 9.22 in \cite{GT} or
theorem 4.8 in \cite{CC}.
Proposition \ref{prop2} is established.

\vskip 5pt
\hfill $\Box$
\vskip 5pt

\noindent{\bf Proof of  Theorem \ref{c11}.}\
We first prove that  $u\in {\cal A}$.  We only need to verify property (A2)
 since property (A1)
 has already been assumed.   Let $\Omega$ be as in the statement of (A2),
 then, by the conformal invariance of $A^u$,
  $A^{ u_{x,\lambda} }$
  and $A^{ (1+\delta)u }$ 
are still in $\partial U$ a.e. in $\Omega$.
Thus, by Proposition \ref{prop2},
(A2) is satisfied.  So we have proved that $u\in {\cal A}$.
By Theorem \ref{thmB15}, (\ref{B15-1}) holds for some $a>0, b\ge 0$ and
$\bar x\in \Bbb R^n$.  We only need to prove that $b=0$.  Suppose that
$b>0$, then $A^u$ is a positive constant multiple of $I$ in $\Bbb R^n$,
and therefore $A^u\in U$ in $\Bbb R^n$ according to  (\ref{A3-0}). 
This violates $A^u\in \partial U$ a.e. in $\Bbb R^n$.
Theorem \ref{c11} is established.

\vskip 5pt
\hfill $\Box$
\vskip 5pt

\section{Proof of Theorem \ref{thmA3new} and Theorem \ref{thm12}}

First we give a variation of Proposition 
\ref{lemA5-2} which allows $u$ to have
isolated singularities.

\begin{prop} Let $U\subset {\cal S}^{n\times n}$ be
an open set satisfying (\ref{C3-1}) and
(\ref{A1-0}),  $F\in C^1(U)\cap C^0(\overline U)$
satisfy (\ref{21-4}) and (\ref{A1-1}), and let 
$\Omega\subset \Bbb R^n$ be a bounded open
set containing $m$ points
$S_m:=\{P_1, \cdots, P_m\}$, $m\ge 1$,
$u\in C^0(\overline \Omega\setminus S_m)$ and
$v\in C^1(\overline \Omega)$ satisfy
\begin{equation}
u>0, \
\Delta u\le 0\ \
\mbox{in}\ \Omega\setminus S_m,\
\
\ v>0\ \mbox{in}\ \overline \Omega,
\label{C2-10}
\end{equation}
\begin{equation}
u>v\quad \mbox{on}\ \partial \Omega.
\label{C2-11}
\end{equation}
We assume that $u$ is a weak solution of
$$
F(A^u)\ge 0\qquad \mbox{in}\
\Omega\setminus S_m,
$$
$v$ is a weak solution of
$$
F(A^v)=0\qquad \mbox{in}\ \Omega\setminus S_m.
$$
Then
\begin{equation}
\inf_{ \Omega \setminus S_m}
(u-v)>0.
\label{C2-12}
\end{equation}
\label{prop4-1}
\end{prop}

\noindent{\bf Proof.}\ We prove it by induction on 
the number of points $m$.  
We start from $m=0$ with 
$S_0=\emptyset$.  The result  is contained in
Proposition \ref{lemA5-2}.  Now we assume
that the result holds for $m-1$ points,
$m-1\ge 0$, and we will prove it for $m$ points.

Let
$$
S_m=\{P_1, \cdots, P_m\}\subset \Omega.
$$
By the assumption on $u$ and $v$, there exist
$\{\beta_i\}$, $\{\tilde \beta_i\}\subset
C^0(\Omega\setminus S_m, {\cal S}^{n\times n})$
and positive functions 
$\{u_i\}$, $\{v_i\}\subset C^2(\Omega\setminus S_m)$ such that,
for any compact subset $K$ of $\Omega\setminus S_m$,
$$
A^{u_i}+\beta_i\in U,\
A^{v_i}+\tilde \beta_i\in U,
\qquad \mbox{in}\ \Omega\setminus S_m,
$$
$$
u_i\to u, \
v_i\to v, \ \beta_i\to 0,\
\tilde \beta_i\to 0,
\qquad \mbox{in}\ C^0(K),
$$
and
$$
F\left(A^{v_i}+\tilde \beta_i\right)\to 0\qquad
\mbox{in}\ C^0(K).
$$

We prove (\ref{C2-12}) by contradiction argument.
Suppose it does not hold, then
$$
\inf_{ \Omega\setminus S_m}(u-v)\le 0.
$$
Since $u>0$ in $\Omega\setminus S_m$,
$\Delta u\le 0$ in $\Omega\setminus S_m$, we know that
$\displaystyle{
\inf_{ \Omega \setminus S_m }u>0.}$
Thus, for some $0<a\le 1$,
$$
\inf_{ \Omega\setminus S_m}(u-av)=0.
$$
Since we can use $a^{-1}u$ instead of $u$, we may assume without loss
of generality that $a=1$.  So we have, in
addition, 
\begin{equation}
\inf_{ \Omega\setminus S_m}(u-v)=0.
\label{C6-10}
\end{equation}

Let $P_m$ be the origin, and let
$$
\widehat \Omega:= \Omega\setminus\{P_1, \cdots, P_{m-1}\}.
$$
We fix some positive function $\varphi$ so that the hypotheses
on $\varphi$ in Theorem \ref{thm10}, with $\widehat \Omega$
being the $\Omega$ there, are satisfied.  Since
$\Phi(v, 0, 1; \cdot)=v$ and
$u>v$ on $\partial \Omega$, we can fix some
small $\epsilon_4>0$ so that
$|x|\le \epsilon_4$ and $|\lambda-1|\le \epsilon_4$
guarentee
\begin{equation}
u>\Phi(v,x,\lambda;\cdot)\qquad
\mbox{on}\ \partial\Omega.
\label{C8-10}
\end{equation}
For such $x$ and $\lambda$, if we assume both
\begin{equation}
\inf_{ \widehat \Omega \setminus \{0\} }
\left[  u-\Phi(v,x,\lambda;\cdot)\right]=0
\label{C8-11}
\end{equation}
and
$$
\liminf_{ |y|\to 0}
\left[ u(y)-\Phi(v,x,\lambda; y)\right]>0,
$$
we would have, for some small $\tilde \epsilon,
\hat \epsilon>0$,
\begin{equation}
u(y)-\Phi(v,x,\lambda; y)>\hat \epsilon>0,
\qquad\forall\ 0<|y|\le \tilde \epsilon.
\label{C9-10}
\end{equation}

Let
$$
\tilde u(y):=\lambda^{ \frac{n-2}2}\varphi(\lambda)^{-1}u(y),
\qquad \tilde v(y):=\lambda^{\frac{n-2}2 }v(x+\lambda y).
$$
We know from (\ref{C8-10}) and (\ref{C9-10}) that
$$
\tilde u>\tilde v\qquad\mbox{on}\
\partial (\Omega \setminus \overline B_{ \tilde \epsilon}).
$$
It is easy to see that the hypotheses of Proposition \ref{prop4-1},
with $u$ replaced by $\tilde u$, $v$ replaced by $\tilde v$,
$\Omega$ replaced by $\Omega \setminus\overline B_{\tilde \epsilon}$,
$S_m$ replaced by $S_{m-1}:=
\{P_1, \cdots, P_{m-1}\}$, are satisfied.  By the induction
hypothesis,
$$
\inf_{  (\Omega \setminus \overline B_{\tilde \epsilon})\setminus S_{m-1}  }
\left(\tilde u-\tilde v\right)>0,
$$
i.e.
$$
\inf_{  (\Omega \setminus \overline B_{\tilde \epsilon})\setminus S_{m-1}  }
\left[u-\Phi(v,x,\lambda; \cdot)\right]>0.
$$
This and (\ref{C9-10}) violate (\ref{C8-11}).
Impossible.  Thus we have proved that (\ref{C8-11})
implies 
\newline
$\displaystyle{
\liminf_{ |y|\to 0}
[u(y)-\Phi(v,x,\lambda;y)]=0}.$
Namely, we have verified (\ref{gg5znew}),
with $\Omega$ replaced by $\widehat \Omega$.  Therefore,
by Theorem \ref{thm10}, either (\ref{3-2y}) holds or
$u=v=v(0)$ near the origin.

If (\ref{3-2y}) holds, then 
$\displaystyle{
\inf_{  B_\epsilon\setminus \{0\} }(u-v)>0 }$
for some $\epsilon>0$.
By the induction hypotheses, applied  on $\Omega\setminus \overline B_\epsilon$,
we obtain
$\displaystyle{
\inf_{  (\Omega\setminus\overline B_\epsilon)\setminus \{P_1, \cdots, P_{m-1}\}  }
(u-v)>0}.$
It follows that $\displaystyle{
\inf_{ \Omega\setminus S_m}(u-v)>0},$
violating (\ref{C6-10}).  A contradiction.

If $u=v=v(0)$ near the origin, say in $B_{\bar \epsilon}$ for
some small $\bar \epsilon>0$, we let
$\displaystyle{
\varphi(y)=e^{ \delta |y|^2 } }$ be 
the function in Lemma \ref{lemv}, and let
$$
u_\epsilon:= \left( u^{ -\frac 2{n-2} }+\epsilon\varphi\right)^{ -\frac {n-2}2 },
\qquad a_\epsilon:=\inf_{  \Omega\setminus \{P_1, \cdots, P_m\} }
\frac{u_\epsilon}v.
$$
Since 
$\displaystyle{
\frac{u_\epsilon}v=\frac uv+O(\epsilon)}$, it
is easy to see from (\ref{C6-10}) that
\begin{equation}
\lim_{ \epsilon\to 0}a_\epsilon=1.
\label{C14-10}
\end{equation}

For $|y|\le \bar \epsilon$, $u(y)=v(y)=v(0)$.
So
$$
\frac{  u_\epsilon(y)  }{  v(y)  }
=1-\frac{n-2}2 \epsilon v(0)^{  \frac 2{n-2} }
e^{ \delta|y|^2 }  +O(\epsilon^2).
$$
For $\epsilon>0$ small,
\begin{equation}
\frac{ u_\epsilon}v >
\max_{  |y|=\bar \epsilon }
\frac{  u_\epsilon(y)  }{  v(y)  }\ge a_\epsilon\qquad
\mbox{on}\
B_{ \bar \epsilon/2 }.
\label{C15-10}
\end{equation}
Taking $\epsilon>0$ smaller if necessary, we have, in 
view of (\ref{C14-10}) and the fact $u>v$ on $\partial \Omega$,
$$
\frac{  u_\epsilon }v > a_\epsilon\qquad
\mbox{on}\ \partial \Omega.
$$
Fix this $\epsilon>0$ and let
$\displaystyle{
\widehat \Omega :=\Omega\setminus \overline B_{\bar \epsilon/2} }$.
We know that
\begin{equation}
(a_\epsilon)^{-1}u_\epsilon>v\qquad
\mbox{on}\ \partial \widehat \Omega.
\label{C16-10}
\end{equation}
Let
$$
u_i^\epsilon:=\left( u_i^{-\frac 2{n-2} }+\epsilon\varphi\right)^{ -\frac {n-2}2 }.
$$
Making $\epsilon>0$ smaller if necessary, we have, by Lemma \ref{lemv},
$\displaystyle{
A^{  (a_\epsilon)^{-1}u_i^\epsilon }\in U }$
in $\widehat \Omega$ for large $i$.
Clearly,
$\displaystyle{
(a_\epsilon)^{-1}u_i^\epsilon  \to
(a_\epsilon)^{-1}u^\epsilon }$
in $C^0_{loc}(\widehat \Omega)$.  By the induction hypothesis,
in view of (\ref{C16-10}), we obtain
$$
\inf_{  \widehat \Omega \setminus
\{P_1, \cdots, P_{m-1}\}  }
\left(  (a_\epsilon)^{-1}u_\epsilon-v\right)>0.
$$
This and (\ref{C15-10}) imply
$$
\inf_{  \Omega\setminus \{P_1, \cdots, P_m\} }
\frac{  u_\epsilon}v >a_\epsilon,
$$
violating the definition of $a_\epsilon$.
Impossible.  Proposition \ref{prop4-1} is established.

  \vskip 5pt
    \hfill $\Box$
     \vskip 5pt

Similar to Proposition \ref{prop4-1}, we have the following variation of
Proposition \ref{prop2}.
\begin{prop} Let $U\subset {\cal S}^{n\times n}$ be
an open set satisfying (\ref{C3-1}) and
(\ref{A1-0}),  $F\in C^1(U)\cap C^0(\overline U)$
satisfy (\ref{21-4}) and (\ref{A1-1}), and let 
$\Omega\subset \Bbb R^n$ be a bounded open
set containing $m$ points
$S_m:=\{P_1, \cdots, P_m\}$, $m\ge 1$,
$u\in C^0(\overline \Omega\setminus S_m)\cap C^{1,1}_{loc}(\Omega\setminus S_m)$ and
$v\in C^1(\overline \Omega)\cap C^{1,1}_{loc}(\Omega\setminus S_m)$ satisfy
(\ref{C2-10}),
(\ref{C2-11}) and 
$$
F(A^u)\ge 0,\ \ 
F(A^v)=0\qquad \mbox{almost everywhere in}\ \Omega\setminus S_m.
$$
Then (\ref{C2-12}) holds.
\label{prop4-2}
\end{prop}

\noindent{\bf Proof.}\ It follows from modification of the proof of
Propsosition \ref{prop4-1}, using Proposition \ref{prop2}
instead of Proposition \ref{lemA5-2}, and making some other obvious changes.  We omit the details.

  \vskip 5pt
    \hfill $\Box$
     \vskip 5pt

Now we give the

\noindent{\bf Proof of Theorem \ref{thmA3new}.}\
By (\ref{A3-1new}) and the positivity of $u$,
$$
\liminf_{|y|\to 0} u(y)>0,
\qquad
\liminf_{|y|\to \infty} |y|^{n-2} u(y)>0.
$$
As usual, for any $x\in \Bbb R^n\setminus\{0\}$ and for any 
$0<\delta<1$,
$$
\bar\lambda_\delta(x)=\sup\{0<\mu<|x|\
|
u_{x,\lambda}(y)\le (1+\delta)u(y),
\forall\ 0<\lambda<\mu, |y-x|>\lambda, |y|\ne 0\}>0
$$
is well defined.

By the definition of 
$\bar\lambda_\delta(x)$, 
\begin{equation}
u_{x, \bar\lambda_\delta(x) }(y)\le (1+\delta)u(y),
\qquad \forall\ |y-x|\ge \bar\lambda_\delta(x), y\ne 0.
\label{211}
\end{equation}

\begin{lem} If $\bar\lambda_\delta(x)<|x|$ for 
some $0<\delta<1$ and $x\in \Bbb R^n\setminus\{0\}$,
then either
$$
\liminf_{ |y|\to 0}
 \left[
(1+\delta)u(y)-u_{x, \bar\lambda_\delta(x)}(y)\right]=0
$$
or
$$
\liminf_{ |y|\to\infty}
|y|^{n-2} \left[
(1+\delta)u(y)-u_{x, \bar\lambda_\delta(x)}(y)\right]=0
$$
\label{lemO4}
\end{lem}

\noindent{\bf Proof.}\ Suppose the contrary, then  for some
$0<\delta<\bar\delta$, $x\in \Bbb R^n\setminus\{0\}$,
$\bar \lambda_\delta(x)<|x|$, and for some
$\displaystyle{
R>3|x|+\frac{ 3 }{|x|-\bar\lambda_\delta(x) }}$, $0<\epsilon_2< \frac 19 \min\{ |x|,
|x|-\bar\lambda_\delta(x)\}$,
$\epsilon>0$,
 $$
 (1+\delta)u(y)-u_{x, \lambda}(y)\ge
 \frac{\epsilon}{|y|^{n-2}},\qquad
 \forall\ |\lambda-\bar\lambda_\delta(x)|<\epsilon_2, |y-x|\ge R,
 $$
$$
 (1+\delta)u(y)-u_{x, \lambda}(y)\ge
 \epsilon, \qquad
  \forall\ |\lambda-\bar\lambda_\delta(x)|<\epsilon_2,
  y\in B_{\frac 1R}\setminus\{0\}\cup \partial B_{ \lambda}.
  $$

Let
$$
\Omega:=\{y\in \Bbb R^n\ |\ \frac 1R<
|y|<R, |y-x|> 
 \bar\lambda_\delta(x)\},
 $$
By Lemma \ref{lemA5-1} and Lemma \ref{lemA10-1}, 
$(1+\delta)u$ is a weak solution of
$$
F\left(A^{(1+\delta)u}\right)=0\qquad \mbox{in}\ \Omega,
$$
and $u_{ x, \bar\lambda(x) }$ is a weak solution
of 
$$
F\left( A^{  u_{x,\bar\lambda_\delta(x)}  }\right)=0\qquad \mbox{in}\ \Omega.
$$
We also know that
$$
(1+\delta)u\ge u_{x,\bar\lambda_\delta(x)}\ \
\mbox{on}\ \overline \Omega \qquad
\mbox{and}\qquad  
(1+\delta)u>
  u_{x,\bar\lambda_\delta(x)}\ \ \mbox{on}\ \partial \Omega.
$$
It follows, using Proposition \ref{lemA5-2}, that
$$
 (1+\delta)u> u_{x, \bar\lambda_\delta(x)}\qquad
 \mbox{on}\ \overline \Omega.
 $$
As usual, the moving sphere procedure can go beyond
$\bar \lambda_\delta(x)$, violating the definition of
$\bar \lambda_\delta(x)$.
Lemma \ref{lemO4} is established.

  \vskip 5pt
    \hfill $\Box$
     \vskip 5pt

\begin{lem} For all $0<\delta<1$ and for all $x\in \Bbb R^n\setminus \{0\}$,
$\bar \lambda_\delta(x)=|x|$.
\label{lem4-2}
\end{lem}

\noindent{\bf Proof.}\ We prove it by contradiction argument.  Suppose
the contrary, then $\bar \lambda_\delta(x)<|x|$ for
some $x\in \Bbb R^n\setminus \{0\}$.  Let
$$
\Omega := B_{ \bar \lambda_\delta(x)  },
\quad P_1:= x-\frac {  \bar\lambda_\delta(x)^2 x }{  |x|^2 },
\quad P_2:= x,
$$
$$
\tilde u(y):= (1+\delta) u_{ x, \bar \lambda_\delta(x) }(y),\qquad
\tilde v(y):= u(y).
$$
Since $0<\bar \lambda_\delta(x)<|x|$, we have $P_1, P_2\in \Omega$.
By (\ref{211}),
$$
\tilde u\ge \tilde v\qquad \mbox{in}\ \Omega\setminus\{P_1, P_2\}.
$$
It is easy to check that the hypotheses of
Proposition \ref{prop4-1}, with $u$ replaced by
$\tilde u$, $v$ replaced by $\tilde v$, and with $m=2$, are satisfied.
Thus, by Proposition \ref{prop4-1}, 
$$
\inf_{  \Omega\setminus\{P_1, P_2\} }(\tilde u-\tilde v)>0.
$$
This implies
$$
\liminf_{ |y|\to 0}
\left[ (1+\delta) u(y) -u_{ x, \bar\lambda_\delta(x) }(y)\right]
=\left(\frac {  \bar\lambda_\delta(x)}{|x|}\right)^{n-2}
\liminf_{z\to P_1}\left[ \tilde u(z)-\tilde v(z)\right]>0,
$$
and
$$
\liminf_{|y|\to \infty} |y|^{n-2}
\left[ (1+\delta) u(y) -u_{ x, \bar\lambda_\delta(x) }(y)\right]
=\left( \bar\lambda_\delta(x)\right)^{n-2}
\liminf_{z\to P_2}\left[ \tilde u(z)-\tilde v(z)\right]>0,
$$
which contradict to Lemma \ref{lemO4}.
Impossible.  Lemma \ref{lem4-2} is established.

  \vskip 5pt
    \hfill $\Box$
     \vskip 5pt

Now we complete the proof of 
 Theorem \ref{thmA3new}.
By Lemma \ref{lem4-2},
$$
\bar\lambda_\delta(x)=|x|,\qquad \forall\ x\in \Bbb R^n\setminus\{0\},
\forall\ 0<\delta<1,
$$
i.e.
$$
u_{x, \lambda}(y)\le (1+\delta)u(y),\qquad
\forall\  0<\lambda<|x|,  |y-x|\ge\lambda, y\ne 0.
$$
Sending $\delta$ to $0$ leads to (\ref{JJJ}). 
The radial symmetry of  $u$ follows from
(\ref{JJJ}), 
 see e.g. \cite{Li1}. 
Theorem \ref{thmA3new} is established.

 \vskip 5pt
     \hfill $\Box$
           \vskip 5pt

\noindent{\bf Proof of Theorem \ref{thm12}.}\
The proof is similar to that of Theorem \ref{thmA3new},
using Proposition \ref{prop2} instead of Proposition \ref{lemA5-2}
and using Proposition \ref{prop4-2} instead of
Proposition \ref{prop4-1}.
We omit the details.

\vskip 5pt
\hfill $\Box$

\end{document}